\documentclass{article}
\usepackage{amsmath}

\usepackage{amssymb}
\usepackage{amsfonts}
\usepackage{color}
\usepackage[bottom]{footmisc}

\input{tcilatex}
\begin{document}

\title{A Class of Orderings in the Range of Borda's Rule}
\author{Jerry S. Kelly\thanks{%
Department of Economics, Syracuse University. E-mail: jskelly@maxwell.syr.edu%
} \and Shaofang Qi\thanks{%
Syracuse University and Humboldt University Berlin (after 1 September 2015).
E-mail: sqi@syr.edu}}
\maketitle

\begin{abstract}
We present a class of orderings $\tciLaplace $ for which there exists a
profile $u$ of preferences for a fixed odd number of individuals such that
Borda's rule maps $u$ to $\tciLaplace $.\bigskip \bigskip \bigskip 
\end{abstract}

Kelly and Qi [1] initiated the study of what orderings are in the range of
Borda's rule for profiles of strong preference orderings and a fixed number
of individuals. Here we extend those results, establishing a new class of
orderings in the Borda range.\medskip 

\textbf{1. Framework.\footnote{%
Most of this section is drawn from Kelly and Qi [1].}}\medskip 

We begin with a finite set $N=\{1,...,n\}$ of individuals, $n\geq 2$, and a
finite set $X$ of \textbf{alternatives}, with $\mid X\mid $ $=m\geq 2$. A 
\textbf{binary relation} $\rho $ on $X$ is a non-empty subset of the
Cartesian product, $X\times X$; if $(x,y)\in \rho $, we will often write $%
x\rho y$. Relation $\rho $ is\medskip

\qquad 1. \textbf{reflexive} if $x\rho x$ for all $x$ in $X$;

\qquad 2. \textbf{asymmetric} if for all $x,y$ in $X$: $x\rho y$ and $y\rho
x $ imply $x=y$;

\qquad 3. \textbf{complete} if for all for all $x,y$ in $X$ such that $x\neq
y$,

\qquad \qquad either $x\rho y$ or $y\rho x$;

\qquad 4. \textbf{transitive} if $x\rho y$ and $y\rho z$ imply $x\rho z$ for
all $x,y,z$ in $X$.\medskip

Relation $\rho $ is a \textbf{weak order} on $X$ if it is a reflexive,
complete, and transitive relation on $X$; $\rho $ is a \textbf{strong order}
on $X$ if it is a weak order on $X$ and is also asymmetric. The set of all
strong orders on $X$ is denoted $L(X)$. If $r$ is a strong order on $X$,
then $r[1]$ is the top-ranked alternative in $r$: $x\succ y$ for all $y$ in $%
X\backslash \{x\}$. More generally, $r[k]$ is the $k^{th}$-ranked
alternative in $r$. \ The \textbf{inverse} $R^{-1}$ of an order $R$ is
defined by $xR^{-1}y$ if and only if $yRx$. \ A \textbf{profile} is an
ordered $n$-tuple $u=(u(1),u(2),...,u(n))\in L(X)^{n}$ of weak
orders.\medskip

Given a profile $u$ in $L(X)^{n}$, define $s(u,x,i)=k$, where $u(i)[k]=x$.\
Then the \textbf{Borda score} of $x$ at $u$, $S(u,x)$, is the sum of the $%
s(u,x,i)$ over $i$, for $1\leq i\leq n$. The \textbf{Borda ranking}, $%
f_{B}(u)$, sets $x\succ y$ if and only if $S(u,x)\leq S(u,y)$, \medskip

The outcome of the Borda ranking procedure is a weak ordering $\tciLaplace $
which could be written as%
\[
\tciLaplace =X_{1}\succ X_{2}\succ ...\succ X_{T}
\]%
\newline
where: (i) each $X_{i}\subset X$, (ii) the $X_{i}$ are pairwise disjoint,
(iii) alternatives within an $X_{i}$ all have the same Borda score, and (iv) 
$i<j$ implies all alternatives in $X_{i}$ have Borda score less than all
alternatives in $X_{j}$. \ Each $X_{i}$ is called a \textbf{level}.\medskip

Because Borda satisfies neutrality, we can usefully abbreviate our
descriptions of weak order images under Borda's rule. \ For a given $X$, we
only have to be concerned with the \textit{number} of alternatives in each
level, not with exactly which alternatives are in each level. Showing $%
\tciLaplace =\{a,b,c\}\succ \{d,e\}\succ \{f\}\succ \{g,h,i,j\}$ is in the
image of $f_{B}$ also shows that $\tciLaplace ^{\ast }=\{i,b,e\}\succ
\{c,h\}\succ \{j\}\succ \{g,e,a,f\}$ is in the range. \ More generally, with
a slight abuse of language, we say a weak order generated by Borda's rule 
\textit{is} a sequence $(m_{1},m_{2},...,m_{T})$ where the $m_{i}$ are the
cardinalities of the sets of alternatives with the same Borda score.\medskip

Kelly and Qi [1] established several propositions regarding the Pareto
range. \ For example, if at least one $m_{i}$ is odd, $%
(m_{1},m_{2},...,m_{T})$ is in the Borda range for all odd $n$.
(Accordingly, if $m$ is odd, $\tciLaplace $ is in the Borda range for all
odd n.) Other results relevant for this paper are:\medskip

\textbf{Theorem 3.} \ Suppose $\tciLaplace =(m_{1},m_{2},...,m_{T})$ has
every $m_{i}$ is even. \ Let $k\geq 1$ be the largest power of $2$ dividing
all the $m_{i}$, so $\tciLaplace =(2^{k}s_{1},,2^{k}s_{2}...,2^{k}s_{T}).$
If $s_{1}+s_{2}+...+s_{T}$ is odd, then for every odd positive integer $n$,
there \underline{does not} exist\ a profile $u$ such that $%
f_{B}(u)=\tciLaplace $.\medskip

\textbf{Lemma 4.} \ Suppose $\tciLaplace =(m_{1},m_{2},...,m_{T})$ has every 
$m_{i}$ is even. \ Let $k\geq 1$ be the largest power of 2 dividing all the $%
m_{i}$, so $\tciLaplace =(2^{k}s_{1},,2^{k}s_{2}...,2^{k}s_{T}).$ If $%
s_{1}+s_{2}+...+s_{T}$ is even and all $s_{i}$ are odd, then for every odd $%
n\geq 3$, there exists\ a profile $u$ such that $f_{B}(u)=\tciLaplace $%
.\medskip

In particular, if $\tciLaplace =(m_{1},m_{2})$ has two equal levels, then $%
\tciLaplace $ is in the Borda range for all odd n.\medskip

\textbf{2. The New Class.}\medskip

Lemma 4 fails to cover most cases where $s_{1}+s_{2}+...+s_{T}$ is even. \
Here we examine long orderings. \ For fixed even $m$, all orders with
sufficiently many levels but with no odd levels will be made up entirely of
levels equal to $2$ or $4$, with not very many $4s$.\medskip

Let $\tciLaplace =(m_{1},m_{2},...,m_{T})$. If all the $m_{i}=2$, then $%
\tciLaplace $ is not in the range of Borda's rule for any odd $n$ if $T$ is
odd (Theorem 3) and is in the range of Borda's rule for all odd $n$ if $T$
is even (Lemma 4). A similar statement can be made if all the $m_{i}=4$%
.\medskip

So we are only interested in the case where \underline{both} $2$ and $4$ do
appear in $\tciLaplace $. By Theorem 3 again, $\tciLaplace $ is not in the
Borda range if there are an odd number of $2s$ in $\tciLaplace $. So we may
suppose $\tciLaplace $ contains $4$ and an even number $(\geq 2)$ of $2s$.
We first treat a set of special cases.\medskip

\qquad \textbf{Lemma}. Each of the following (patterns of) orders (with
exactly two levels equal to $2$) is in the Borda range for $n=3$ (and so for
all odd $n\geq 3$):\medskip

\qquad \qquad \qquad \qquad $(2,4,4,...,4,4,2)$

\qquad \qquad \qquad \qquad $(4,2,4,4,...,4,4,2)$

\qquad \qquad \qquad \qquad $(2,4,4,...,4,4,2,4)$

\qquad \qquad \qquad \qquad $(4,2,4,4,...,4,4,2,4)$\medskip

\qquad A proof of the Lemma appears in the next section.\medskip

\qquad \textbf{Theorem.} Suppose $\tciLaplace =(m_{1},m_{2},...,m_{T})$ and $%
\tciLaplace $ contains only $4s$ and an even number $(\geq 2)$ of $2s$. Then 
$\tciLaplace $ is in the Borda range for all odd $n\geq 3$.\medskip

\qquad \textbf{Proof of the Theorem:} The proof is by induction on the
number of levels $T$ $(\geq 2)$.\medskip

\qquad \textbf{Basis:} For $T=2$, $\tciLaplace $ must be $(2,2)$ and this is
in the range by Lemma 4.\medskip

\qquad \textbf{Induction step:} We now assume the result is true for all
levels less than $T\geq 3$. We first decompose $\tciLaplace $ by stripping
out some $4s$ that might occur at the beginning. Let $\tciLaplace _{0}$ be
the order made up of the largest even sequence of $4s$ prior to the first
occurrence of 2 in $\tciLaplace $. So $\tciLaplace =\tciLaplace _{0}\succ
\tciLaplace _{1}$ where $\tciLaplace _{0}$ contains an even number $(\geq 2)$
of $4s$ and order $\tciLaplace _{1}$ looks like either $(2,...)$ or $%
(4,2,...)$. By Lemma 4, $\tciLaplace _{0}$ is in the range and we will
catenate the profile for $\tciLaplace _{0}$ with the profile we will
construct for $\tciLaplace _{1}$. Let $\tciLaplace _{1}^{\ast }$ be the
initial sequence of $\tciLaplace _{1}$ up to and including the second
occurrence of 2: $(2,4,4,...,4,4,2)$ or $(4,2,4,4,...,4,4,2)$. Now $%
\tciLaplace _{1}=\tciLaplace _{1}^{\ast }\succ \tciLaplace _{2}$ , where $%
\tciLaplace _{2}$ contains an even number of only $2s$ (possibly 0). If $%
\tciLaplace _{2}$ contains a positive even number of $2s$, the induction
hypothesis shows $\tciLaplace _{2}$ is in the Borda range. By the Lemma, $%
\tciLaplace _{1}^{\ast }$ is in the Borda range and we can catenate profiles
to show that $\tciLaplace =\tciLaplace _{0}\succ \tciLaplace _{1}^{\ast
}\succ \tciLaplace _{2}$ is in the Borda range.\medskip

\qquad So suppose $\tciLaplace _{2}$ doesn't contain any $2s$. Either it is
empty and $\tciLaplace =\tciLaplace _{0}\succ \tciLaplace _{1}^{\ast }$ is
in the range by catenation, or $\tciLaplace _{2}=(4,4,...,4)$. If $%
\tciLaplace _{2}$ contains an even number of $4s$, then $\tciLaplace _{2}$
is in the range by Lemma 4 and $\tciLaplace =\tciLaplace _{0}\succ
\tciLaplace _{1}^{\ast }\succ \tciLaplace _{2}$ is in the Borda range by
catenation.\medskip

\qquad All that remains is the case where $\tciLaplace _{2}$ contains an odd
number of 4s. In that case, move one $4$ from $\tciLaplace _{2}$ to $%
\tciLaplace _{1}^{\ast }$ so that now we take $\tciLaplace _{1}^{\ast
}=(2,4,4,...,4,2,4)$ or $(4,2,4,4,...,4,2,4)$ (either order in the Borda
range by the Lemma) and $\tciLaplace _{2}$ has an even number of $4s$. Now, $%
\tciLaplace =\tciLaplace _{0}\succ \tciLaplace _{1}^{\ast }\succ \tciLaplace
_{2}$ where each part is in the Borda range and catenation yields our
result. \ \ \ \ \ \U{25a1}\medskip \medskip 

\textbf{3. \ Proof of the Lemma.}\medskip

We wish to show that each of the following (patterns of) orders (with
exactly two levels equal to $2$) is in the Borda range for $n=3$ (and so for
all odd $n\geq 3$):\medskip

\qquad \qquad \qquad \qquad $(2,4,4,...,4,4,2)$

\qquad \qquad \qquad \qquad $(4,2,4,4,...,4,4,2)$

\qquad \qquad \qquad \qquad $(2,4,4,...,4,4,2,4)$

\qquad \qquad \qquad \qquad $(4,2,4,4,...,4,4,2,4)$\medskip

To construct the relevant profile for each sequence, we will make use of the
profile below (the profile we constructed for Lemma 4 in the main paper). In
particular, for any weak ordering $\tciLaplace =(2s_{1},,2s_{2})$, where
both $s_{1}$ and $s_{2}$ are odd numbers, the following profile $v$ has $%
f_{B}(v)=\tciLaplace $. Note that at profile $v$, for individual \#2, below
the second group of blue options are all odd-subscript options; for
individual \#3, below the second group of blue options are all
even-subscript options. All the blue options (all together $2s_{2}$)
represent the options that have equal Borda score, which is smaller than the
equal Borda score of the remaining options (all together $2s_{1}$). The
Borda score difference between the two groups of options is $\frac{%
s_{1}+s_{2}}{2}$. Individual \#1 ranks options in order from $x_{1}$ to $%
x_{2(s_{1}+s_{2})}$. The profile we will construct for each sequence is
based on variations of $v$. The variations will only be made regarding
individual \#1's ranking of options, and we retain the ranking of individual
\#2 and \#3 unchanged. Therefore, for the following analysis, for
simplicity, we will only re-state individual \#1's ranking at a profile
without presenting the full profile.\medskip 
\[
\begin{tabular}{|lll|}
\hline
$1$ & $2$ & $3$ \\ \hline
$x_{1}$ & \multicolumn{1}{|l}{$x_{(s_{1}+s_{2})+s_{1}-1}$} & $%
x_{(s_{1}+s_{2})+s_{1}}$ \\ 
$x_{2}$ & \multicolumn{1}{|l}{$\vdots $} & $x_{(s_{1}+s_{2})+s_{1}-2}$ \\ 
$\vdots $ & \multicolumn{1}{|l}{$\vdots $} & $\vdots $ \\ 
$\vdots $ & \multicolumn{1}{|l}{$x_{(s_{1}+s_{2})+4}$} & $\vdots $ \\ 
& \multicolumn{1}{|l}{$x_{(s_{1}+s_{2})+2}$} & $x_{(s_{1}+s_{2})+3}$ \\ 
$\vdots $ & \multicolumn{1}{|l}{$%
\color{blue}%
x_{2(s_{1}+s_{2})}$} & $x_{(s_{1}+s_{2})+1}$ \\ 
$\vdots $ & \multicolumn{1}{|l}{$%
\color{blue}%
x_{2(s_{1}+s_{2})-2}$} & $%
\color{blue}%
x_{2(s_{1}+s_{2})-1}$ \\ 
& \multicolumn{1}{|l}{$\vdots $} & $%
\color{blue}%
x_{2(s_{1}+s_{2})-3}$ \\ 
$\vdots $ & \multicolumn{1}{|l}{$\vdots $} & $\vdots $ \\ 
$\vdots $ & \multicolumn{1}{|l}{$%
\color{blue}%
x_{(s_{1}+s_{2})+s_{1}+3}$} & $\vdots $ \\ 
& \multicolumn{1}{|l}{$%
\color{blue}%
x_{(s_{1}+s_{2})+s_{1}+1}$} & $%
\color{blue}%
x_{(s_{1}+s_{2})+s_{1}+2}$ \\ 
$\vdots $ & \multicolumn{1}{|l}{$x_{s_{1}-1}$} & $x_{s_{1}}$ \\ 
$\vdots $ & \multicolumn{1}{|l}{$\vdots $} & $x_{s_{1}-2}$ \\ 
& \multicolumn{1}{|l}{$\vdots $} & $\vdots $ \\ 
$\vdots $ & \multicolumn{1}{|l}{$x_{4}$} & $\vdots $ \\ 
$\vdots $ & \multicolumn{1}{|l}{$x_{2}$} & $x_{3}$ \\ 
& \multicolumn{1}{|l}{$%
\color{blue}%
x_{(s_{1}+s_{2})}$} & $x_{1}$ \\ 
$\vdots $ & \multicolumn{1}{|l}{$%
\color{blue}%
x_{(s_{1}+s_{2})-2}$} & $%
\color{blue}%
x_{(s_{1}+s_{2})-1}$ \\ 
$\vdots $ & \multicolumn{1}{|l}{$\vdots $} & $%
\color{blue}%
x_{(s_{1}+s_{2})-3}$ \\ 
& \multicolumn{1}{|l}{$\vdots $} & $\vdots $ \\ 
$\vdots $ & \multicolumn{1}{|l}{$%
\color{blue}%
x_{s_{1}+3}$} & $\vdots $ \\ 
$\vdots $ & \multicolumn{1}{|l}{$%
\color{blue}%
x_{s_{1}+1}$} & $%
\color{blue}%
x_{s_{1}+2}$ \\ 
& \multicolumn{1}{|l}{$x_{2(s_{1}+s_{2})-1}$} & $x_{2(s_{1}+s_{2})}$ \\ 
$\vdots $ & \multicolumn{1}{|l}{$x_{2(s_{1}+s_{2})-3}$} & $%
x_{2(s_{1}+s_{2})-2}$ \\ 
$\vdots $ & \multicolumn{1}{|l}{$\vdots $} & $\vdots $ \\ 
& \multicolumn{1}{|l}{$\vdots $} & $\vdots $ \\ 
$x_{2(s_{1}+s_{2})-1}$ & \multicolumn{1}{|l}{$x_{3}$} & $x_{4}$ \\ 
$x_{2(s_{1}+s_{2})}$ & \multicolumn{1}{|l}{$x_{1}$} & $x_{2}$ \\ \hline
\end{tabular}%
\]%
\pagebreak

\textbf{I. The sequence }$(2,4,4,...,4,4,2)$.\medskip

\textbf{Case 1.} There are $2k$ (i.e., even number of) $4s$ between the $2s$%
.\medskip

We first construct, using Lemma 4's method, a profile $v$ for the two-level
ordering $(2+\underset{k}{\underbrace{4+\cdots +4}},\underset{k}{\underbrace{%
4+\cdots +4}}+2)$, which is $(2+4k,4k+2)$, which also is $(2(1+2k),2(2k+1))$%
. Given the proof of Lemma 4, the first level consists of two groups of
options, each having $(1+2k)$ options: 
\begin{eqnarray*}
&&\{x_{1},x_{2},\dots ,x_{2k+1}\} \\
&&\{x_{(4k+2)+1},x_{(4k+2)+2},\dots ,x_{(4k+2)+2k+1}\}
\end{eqnarray*}

Similarly, the second level also consists of two groups of options, each
having $(1+2k)$ options: 
\begin{eqnarray*}
&&\{x_{(2k+1)+1},x_{(2k+1)+2},\dots ,x_{(2k+1)+2k+1}\} \\
&&\{x_{(6k+3)+1},x_{(6k+3)+2},\dots ,x_{(6k+3)+2k+1}\}
\end{eqnarray*}

In addition, as Lemma 4 in Kelly and Qi (2015) shows, at $v$, individual
\#1's ranking of options is simply:\medskip 
\[
\begin{tabular}{|l|}
\hline
$1$ \\ \hline
$x_{1}$ \\ 
$x_{2}$ \\ 
$\vdots $ \\ 
$\vdots $ \\ 
$x_{8k+3}$ \\ 
$x_{8k+4}$ \\ \hline
\end{tabular}%
\]

We color individual \#1's above ranking such that options in the first level
are in black and options in the second level are in blue:\medskip 
\[
\begin{tabular}{|l|}
\hline
$1$ \\ \hline
$x_{1}$ \\ 
$x_{2}$ \\ 
$\vdots $ \\ 
$x_{2k+1}$ \\ 
$%
\color{blue}%
x_{(2k+1)+1}$ \\ 
$%
\color{blue}%
x_{(2k+1)+2}$ \\ 
$%
\color{blue}%
\vdots $ \\ 
$%
\color{blue}%
x_{(2k+1)+2k+1}$ \\ 
$x_{(4k+2)+1}$ \\ 
$x_{(4k+2)+2}$ \\ 
$\vdots $ \\ 
$x_{(4k+2)+2k+1}$ \\ 
$%
\color{blue}%
x_{(6k+3)+1}$ \\ 
$%
\color{blue}%
x_{(6k+3)+2}$ \\ 
$%
\color{blue}%
\vdots $ \\ 
$%
\color{blue}%
x_{(6k+3)+2k+1}$ \\ \hline
\end{tabular}%
\]

Now we make changes for (only) individual \#1's ranking to obtain profile $u$
for the level $(2,4,4,...,4,4,2)$ (where there are $2k$ occurrences of $4$).
For the options that are colored the same, we will make the same change. So
we take the first half options, a group of black and a group of blue, as an
example:\medskip 
\[
\begin{tabular}{|l|}
\hline
$1$ \\ \hline
$x_{1}$ \\ 
$x_{2}$ \\ 
$\vdots $ \\ 
$x_{2k+1}$ \\ 
$%
\color{blue}%
x_{(2k+1)+1}$ \\ 
$%
\color{blue}%
x_{(2k+1)+2}$ \\ 
$%
\color{blue}%
\vdots $ \\ 
$%
\color{blue}%
x_{(2k+1)+2k+1}$ \\ 
$\vdots $ \\ \hline
\end{tabular}%
\]

For the options $\{x_{1},x_{2},\dots ,x_{2k+1}\}$, we change the initial
ranking of $x_{1}\succ x_{2}\succ \cdots \succ x_{2k+1}$ to:\medskip 
\[
\begin{tabular}{|l|}
\hline
$%
\color{red}%
x_{2k+1}$ \\ 
$x_{2k-1}$ \\ 
$x_{2k}$ \\ 
$\vdots $ \\ 
$x_{3}$ \\ 
$x_{4}$ \\ 
$x_{1}$ \\ 
$x_{2}$ \\ \hline
\end{tabular}%
\]

And it is straightforward to check that after making the above change, the
options are split into $(k+1)$ levels:\medskip 
\begin{eqnarray*}
&&\{x_{2k+1}\} \\
&&\{x_{2k-1},x_{2k}\} \\
&&\{x_{2k-3},x_{2k-2}\} \\
&&\vdots \\
&&\{x_{3},x_{4}\} \\
&&\{x_{1},x_{2}\}
\end{eqnarray*}%
\medskip

For the options $\{x_{(2k+1)+1},x_{(2k+1)+2},\dots ,x_{(2k+1)+2k+1}\}$, we
change the initial ranking of $x_{(2k+1)+1}\succ x_{(2k+1)+2}\succ \cdots
\succ x_{(2k+1)+2k+1}$ to:\medskip 
\[
\begin{tabular}{|l|}
\hline
$%
\color{blue}%
x_{(2k+1)+2k}$ \\ 
$%
\color{blue}%
x_{(2k+1)+2k+1}$ \\ 
$%
\color{blue}%
x_{(2k+1)+2k-2}$ \\ 
$%
\color{blue}%
x_{(2k+1)+2k-1}$ \\ 
$%
\color{blue}%
\vdots $ \\ 
$%
\color{blue}%
x_{(2k+1)+2}$ \\ 
$%
\color{blue}%
x_{(2k+1)+3}$ \\ 
$%
\color{red}%
x_{(2k+1)+1}$ \\ \hline
\end{tabular}%
\]%
\medskip

And it is straightforward to check that after making the above change, the
options are split into $(k+1)$ levels:\medskip 
\begin{eqnarray*}
&&\{%
\color{blue}%
x_{(2k+1)+2k},%
\color{blue}%
x_{(2k+1)+2k+1}%
\color{black}%
\} \\
&&\{%
\color{blue}%
x_{(2k+1)+2k-2},%
\color{blue}%
x_{(2k+1)+2k-1}%
\color{black}%
\} \\
&&\vdots \\
&&\{%
\color{blue}%
x_{(2k+1)+4},%
\color{blue}%
x_{(2k+1)+5}%
\color{black}%
\} \\
&&\{%
\color{blue}%
x_{(2k+1)+2},%
\color{blue}%
x_{(2k+1)+3}%
\color{black}%
\} \\
&&\{%
\color{blue}%
x_{(2k+1)+1}%
\color{black}%
\}
\end{eqnarray*}

Now we compare the earlier $(k+1)$ levels consisting of $\{x_{1},x_{2},\dots
,x_{2k+1}\}$ with the above $(k+1)$ levels consisting of $%
\{x_{(2k+1)+1},x_{(2k+1)+2},\dots ,x_{(2k+1)+2k+1}\}$.\medskip

Compared with $v$ (before making change of \#1's ranking), for the option $%
x_{2k+1}$ at the highest level of the $(k+1)$ levels consisting of $%
\{x_{1},x_{2},\dots ,x_{2k+1}\}$,\ the Borda score of $x_{2k+1}$ is
increased by $2k$. Compared with $v$, for the options $x_{(2k+1)+2k}$ and $%
x_{(2k+1)+2k+1}$\ at the highest level of the $(k+1)$ levels consisting of $%
\{x_{(2k+1)+1},x_{(2k+1)+2},\dots ,x_{(2k+1)+2k+1}\}$,\ the Borda score of $%
x_{(2k+1)+2k}$ and $x_{(2k+1)+2k+1}$ is increased by $2k-1$. In addition, it
follows from Lemma 4 that at $v$, the Borda score of $x_{1},x_{2},\dots
,x_{2k+1}$ (black) is larger than that of $x_{(2k+1)+1},x_{(2k+1)+2},\dots
,x_{(2k+1)+2k+1}$ (blue) by $(2k+1)$. So the option $x_{2k+1}$ is at the
highest level among the $2(k+1)$ levels.\medskip

Similarly, compared with $v$, for the option $x_{(2k+1)+1}$ at the lowest
level of the $(k+1)$ levels consisting of $\{x_{(2k+1)+1},x_{(2k+1)+2},\dots
,x_{(2k+1)+2k+1}\}$,\ the Borda score of $x_{(2k+1)+1}$ is decreased by $2k$%
. Compared with $v$, for the options $x_{1}$ and $x_{2}$ at the lowest level
of the $(k+1)$ levels consisting of $\{x_{1},x_{2},\dots ,x_{2k+1}\}$,\ the
Borda score of $x_{1}$ and $x_{2}$ is decreased by $2k-1$. So the option $%
x_{(2k+1)+1}$ is at the lowest level among the $2(k+1)$ levels.\medskip

Except $x_{2k+1}$ at the highest level and $x_{(2k+1)+1}$ at the lowest
level, we need to show that the remaining options are actually split into $%
2k $ levels. We do this by calculating the Borda score change from $v$ to $u$%
. For the $k$ levels consisting of $\{x_{1},x_{2},\dots ,x_{2k}\}$:\medskip 
\begin{eqnarray*}
\{x_{2k-1},x_{2k}\} &:&(2k-3)\,more \\
\{x_{2k-3},x_{2k-2}\} &:&(2k-7)\,more \\
&&\vdots \\
\{x_{3},x_{4}\} &:&(5-2k)\,more\,[i.e.,(2k-5)\,less] \\
\{x_{1},x_{2}\} &:&(1-2k)\,more\,[i.e.,(2k-1)\,less]
\end{eqnarray*}

For the $k$ levels consisting of $\{x_{(2k+1)+2},\dots ,x_{(2k+1)+2k+1}\}$%
:\medskip 
\begin{eqnarray*}
\{%
\color{blue}%
x_{(2k+1)+2k},%
\color{blue}%
x_{(2k+1)+2k+1}%
\color{black}%
\} &:&(2k-1)\,more \\
\{%
\color{blue}%
x_{(2k+1)+2k-2},%
\color{blue}%
x_{(2k+1)+2k-1}%
\color{black}%
\} &:&(2k-5)\,more \\
&&\vdots \\
\{%
\color{blue}%
x_{(2k+1)+4},%
\color{blue}%
x_{(2k+1)+5}%
\color{black}%
\} &:&(7-2k)\,more\,[i.e.,(2k-7)\,less] \\
\{%
\color{blue}%
x_{(2k+1)+2},%
\color{blue}%
x_{(2k+1)+3}%
\color{black}%
\} &:&(3-2k)\,more\,[i.e.,(2k-3)\,less]
\end{eqnarray*}

And at $v$, the black options Borda score is larger than the of blue by $%
(2k+1)$. It follows that the Borda score difference between the black
options based on ranking at $u$ and the blue options based on ranking at $v$
is an even number, and since the Borda score change of blue options from $v$
to $u$ is an odd number, no level among the $k$ consisting of $%
\{x_{1},x_{2},\dots ,x_{2k}\}$ will have the same Borda score as any level
among the $k$ consisting of $\{x_{1},x_{2},\dots ,x_{2k}\}$. Therefore, they
remain in $2k$ levels.\medskip

Combining with the fact that $x_{2k+1}$ at the highest level and $%
x_{(2k+1)+1}$ at the lowest level, at $u$, for these $(4k+2)$ options, we
obtain a pattern $(1,2,2,...,2,2,1)$ (where there are $2k$ number of $2$s).
Recall that we also do the same thing for another half of the options at $v$
for individual $\#1$,\newline
$\{x_{(4k+2)+1},\ldots ,x_{(4k+2)+2k+1},%
\color{blue}%
x_{(6k+3)+1},\dots ,%
\color{blue}%
x_{(6k+3)+2k+1}%
\color{black}%
\}$, so actually, at $u$, we have $f_{B}(u)=(2,4,4,...,4,4,2)$ (where there
are $2k$ occurrences of $4$).\medskip

\textbf{Case 2.} There are $2k+1$ (i.e., odd number of) $4s$ between the $2s$%
.\medskip

We follow the same procedure as Case 1, except that we make the following
revisions.\medskip

First, at the initial construction, profile $v$ is for the ordering with two
levels: $(2+\underset{k+1}{\underbrace{4+\cdots +4}},\underset{k}{%
\underbrace{4+\cdots +4}}+2)$, which is $(2+4(k+1),4k+2)$, which also is $%
(2(3+2k),2(2k+1))$. So that according to Lemma 4's method, the first level
consists of two groups of options, each having $(3+2k)$ options: 
\begin{eqnarray*}
&&\{x_{1},x_{2},\dots ,x_{2k+1},x_{2k+2},x_{2k+3}\} \\
&&\{x_{(4k+4)+1},x_{(4k+4)+2},\dots
,x_{(4k+4)+2k+1},x_{(4k+4)+2k+2},x_{(4k+4)+2k+3}\}
\end{eqnarray*}

The second level also consists of two groups of options, each having $(1+2k)$
options: 
\begin{eqnarray*}
&&\{%
\color{blue}%
x_{(2k+3)+1},x_{(2k+3)+2},\dots ,x_{(2k+3)+2k+1}%
\color{black}%
\} \\
&&\{%
\color{blue}%
x_{(6k+7)+1},x_{(6k+7)+2},\dots ,x_{(6k+7)+2k+1}%
\color{black}%
\}
\end{eqnarray*}

Following Case 1, we make changes for individual \#1's ranking at $v$. The
same change will be made for the options in the same level (color), so we
take the first half of options, one group of black option and one group of
blue option, as an example:\medskip 
\[
\begin{tabular}{|l|}
\hline
$1$ \\ \hline
$x_{1}$ \\ 
$x_{2}$ \\ 
$\vdots $ \\ 
$x_{2k+2}$ \\ 
$x_{2k+3}$ \\ 
$%
\color{blue}%
x_{(2k+3)+1}$ \\ 
$%
\color{blue}%
x_{(2k+3)+2}$ \\ 
$%
\color{blue}%
\vdots $ \\ 
$%
\color{blue}%
x_{(2k+3)+2k+1}$ \\ 
$\vdots $ \\ \hline
\end{tabular}%
\]

For the options $\{x_{1},x_{2},\dots ,x_{2k+3}\}$, we change the initial
ranking of $x_{1}\succ x_{2}\succ \cdots \succ x_{2k+3}$ of individual \#1
at $v$ to:\medskip 
\[
\begin{tabular}{|l|}
\hline
$%
\color{red}%
x_{2k+3}$ \\ 
$x_{2k+1}$ \\ 
$x_{2k+2}$ \\ 
$\vdots $ \\ 
$x_{3}$ \\ 
$x_{4}$ \\ 
$x_{1}$ \\ 
$x_{2}$ \\ \hline
\end{tabular}%
\]

And the options are split into $(k+2)$ levels, with the Borda score change
going from $v$ to $u$ marked as below:\medskip 
\begin{eqnarray*}
\{x_{2k+3}\} &:&(2k+2)\,more \\
\{x_{2k+1},x_{2k+2}\} &:&(2k-1)\,more \\
\{x_{2k-1},x_{2k}\} &:&(2k-5)\,more \\
&&\vdots \\
\{x_{3},x_{4}\} &:&(3-2k)\,more\,[i.e.,(2k-3)\,less] \\
\{x_{1},x_{2}\} &:&(-1-2k)\,more\,[i.e.,(2k+1)\,less]
\end{eqnarray*}

Similarly, for the options $\{x_{(2k+3)+1},x_{(2k+3)+2},\dots
,x_{(2k+3)+2k+1}\}$, we change the initial ranking of $x_{(2k+3)+1}\succ
x_{(2k+3)+2}\succ \cdots \succ x_{(2k+3)+2k+1}$ to:\medskip 
\[
\begin{tabular}{|l|}
\hline
$%
\color{blue}%
x_{(2k+3)+2k}$ \\ 
$%
\color{blue}%
x_{(2k+3)+2k+1}$ \\ 
$%
\color{blue}%
x_{(2k+3)+2k-2}$ \\ 
$%
\color{blue}%
x_{(2k+3)+2k-1}$ \\ 
$%
\color{blue}%
\vdots $ \\ 
$%
\color{blue}%
x_{(2k+3)+2}$ \\ 
$%
\color{blue}%
x_{(2k+3)+3}$ \\ 
$%
\color{red}%
x_{(2k+3)+1}$ \\ \hline
\end{tabular}%
\]

And the options are split into $(k+1)$ levels, with the Borda score change
going from $v$ to $u$ marked as below:\medskip 
\begin{eqnarray*}
\{%
\color{blue}%
x_{(2k+3)+2k},%
\color{blue}%
x_{(2k+3)+2k+1}%
\color{black}%
\} &:&(2k-1)\,more \\
\{%
\color{blue}%
x_{(2k+3)+2k-2},%
\color{blue}%
x_{(2k+3)+2k-1}%
\color{black}%
\} &:&(2k-5)\,more \\
&&\vdots \\
\{%
\color{blue}%
x_{(2k+3)+4},%
\color{blue}%
x_{(2k+3)+5}%
\color{black}%
\} &:&(7-2k)\,more\,[i.e.,(2k-7)\,less] \\
\{%
\color{blue}%
x_{(2k+3)+2},%
\color{blue}%
x_{(2k+3)+3}%
\color{black}%
\} &:&(3-2k)\,more\,[i.e.,(2k-3)\,less] \\
\{%
\color{blue}%
x_{(2k+3)+1}%
\color{black}%
\} &:&-2k\,more\,[i.e.,2k\,less]
\end{eqnarray*}

Now we compare, at $u$, the $(k+2)$-level group of options (black) with the $%
(k+1)$-level group of options (blue).\medskip

Recall that at $v$, where the black options are a single level and the blue
options are a single level, the Borda score difference (the black ones have
larger Borda score) between the two levels is $\frac{(3+2k)+(2k+1)}{2}=2k+2$%
.\medskip

It follows that option $x_{2k+3}$ is at the highest level among the $(2k+3)$
levels, and option $x_{(2k+3)+1}$ is at the lowest level among the $(2k+3)$
levels.\medskip

To compare the remaining options, we summarize the Borda score difference
between the black options based on ranking at $u$ and the blue options based
on ranking at $v$ (therefore a single level for the blue options)
below:\medskip 
\begin{eqnarray*}
\{x_{2k+3}\} &:&(2k+2)+(2k+2)\,more\,[i.e.,(4k+4)\,more] \\
\{x_{2k+1},x_{2k+2}\} &:&(2k-1)+(2k+2)\,more\,[i.e.,(4k+1)\,more] \\
\{x_{2k-1},x_{2k}\} &:&(2k-5)+(2k+2)\,more\,[i.e.,(4k-3)\,more] \\
&&\vdots \\
\{x_{3},x_{4}\} &:&(3-2k)+(2k+2)\,more\,[i.e.,5\,more] \\
\{x_{1},x_{2}\} &:&(-1-2k)+(2k+2)\,more\,[i.e.,1\,more]
\end{eqnarray*}

We pick the level among the $(k+1)$-level group of blue options where the
Borda score change for the options at the level from $v$ to $u$ is equal to $%
1$ (if no such level exists, then we skip this step). For example, when $k=3$%
, such level consists of options $\{%
\color{blue}%
x_{13},%
\color{blue}%
x_{14}%
\color{black}%
\}$. We then move $x_{1}$ and $x_{2}$ right below the options for this
level. If after this change, the Borda score of $x_{1}$ and $x_{2}$ is still
equal to some level among the $(k+1)$-level group of blue options, then move 
$x_{1}$ and $x_{2}$ further down by two more options. Accordingly, one can
check that the obtained profile has $(2k+1)$ levels among these
options.\medskip

Again, we do the same changes for the second half of the options for
individual \#1. For the obtained profile, the Borda score ordering is $%
(2,4,4,...,4,4,2)$ (where there are $2k+1$ occurrences of $4$).\bigskip

\textbf{II. The sequence }$(4,2,4,...,4,4,2)$.\medskip

\textbf{Case 1.} There are $2k$ (i.e., even number of) occurrences of $4$%
.\medskip

We follow the same steps as in case 1 in the first sequence, until for
individual \#1 at profile $v$, for the options $\{x_{1},x_{2},\dots
,x_{2k+1}\}$, we change the initial ranking of $x_{1}\succ x_{2}\succ \cdots
\succ x_{2k+1}$ to:\medskip 
\[
\begin{tabular}{|l|}
\hline
$%
\color{red}%
x_{2k}$ \\ 
$%
\color{red}%
x_{2k+1}$ \\ 
$%
\color{green}%
x_{2k-1}$ \\ 
$x_{2k-3}$ \\ 
$x_{2k-2}$ \\ 
$\vdots $ \\ 
$x_{3}$ \\ 
$x_{4}$ \\ 
$x_{1}$ \\ 
$x_{2}$ \\ \hline
\end{tabular}%
\]

Here the options are split into $(k+1)$ levels, with the Borda score change
from $v$ to $u$ marked as below:\medskip 
\begin{eqnarray*}
\{x_{2k},x_{2k+1}\} &:&(2k-1)\,more \\
\{x_{2k-1}\} &:&(2k-4)\,more \\
\{x_{2k-3},x_{2k-2}\} &:&(2k-7)\,more \\
\{x_{2k-5},x_{2k-4}\} &:&(2k-11)\,more \\
&&\vdots \\
\{x_{3},x_{4}\} &:&(5-2k)\,more\,[i.e.,(2k-5)\,less] \\
\{x_{1},x_{2}\} &:&(1-2k)\,more\,[i.e.,(2k-1)\,less]
\end{eqnarray*}

Recall that from Lemma 4, at $v$, the Borda score of $x_{1},x_{2},\dots
,x_{2k+1}$ (black) is larger than that of $x_{(2k+1)+1},x_{(2k+1)+2},\dots
,x_{(2k+1)+2k+1}$ (blue) by $(2k+1)$. So the Borda score difference between
the black options based on ranking at $u$ and the blue options based on
ranking at $v$ (therefore a single level for the blue options) is:\medskip 
\begin{eqnarray*}
\{x_{2k},x_{2k+1}\} &:&(2k-1)+(2k+1)\,more\,[i.e.,4k\,more] \\
\{x_{2k-1}\} &:&(2k-4)+(2k+1)\,more\,[i.e.,(4k-3)\,more] \\
\{x_{2k-3},x_{2k-2}\} &:&(2k-7)+(2k+1)\,more\,[i.e.,(4k-6)\,more] \\
\{x_{2k-5},x_{2k-4}\} &:&(2k-11)+(2k+1)\,more\,[i.e.,(4k-10)\,more] \\
&&\vdots \\
\{x_{3},x_{4}\} &:&(5-2k)+(2k+1)\,more\,[i.e.,6\,more] \\
\{x_{1},x_{2}\} &:&(1-2k)+(2k+1)\,more\,[i.e.,2\,more]
\end{eqnarray*}

For options $\{x_{(2k+1)+1},\dots ,x_{(2k+1)+2k+1}\}$, we make the exactly
same change as for case 1 of sequence I, so the options are split into $%
(k+1) $ and the Borda score change going from $v$ to $u$ is:\medskip 
\begin{eqnarray*}
\{%
\color{blue}%
x_{(2k+1)+2k},%
\color{blue}%
x_{(2k+1)+2k+1}%
\color{black}%
\} &:&(2k-1)\,more \\
\{%
\color{blue}%
x_{(2k+1)+2k-2},%
\color{blue}%
x_{(2k+1)+2k-1}%
\color{black}%
\} &:&(2k-5)\,more \\
&&\vdots \\
\{%
\color{blue}%
x_{(2k+1)+4},%
\color{blue}%
x_{(2k+1)+5}%
\color{black}%
\} &:&(7-2k)\,more\,[i.e.,(2k-7)\,less] \\
\{%
\color{blue}%
x_{(2k+1)+2},%
\color{blue}%
x_{(2k+1)+3}%
\color{black}%
\} &:&(3-2k)\,more\,[i.e.,(2k-3)\,less] \\
\{%
\color{blue}%
x_{(2k+1)+1}%
\color{black}%
\} &:&-2k\,more\,[i.e.,2k\,less]
\end{eqnarray*}

We then compare the above two $(k+1)$-level groups of options.\medskip

For $k=1$, the sequence is $(4,2,4,2)$, which is shown to be in the Borda
range for $n=3$ in the Appendix.\medskip

For $k>1$, we have $(4k-3)>(2k-1)$, so that the two levels, $%
\{x_{2k},x_{2k+1}\}$ and $\{x_{2k-1}\}$, remain as highest levels among the $%
2(k+1)$ levels. Similarly, the level $\{%
\color{blue}%
x_{(2k+1)+1}%
\color{black}%
\}$ remains as the lowest one among the $2(k+1)$ levels. For the remaining $%
(2k-1)$ levels, the Borda score change for the black options is even while
for blue options is odd, so they remain as $(2k-1)$ levels. By making the
same change for the second half of the options for individual \#1, we obtain
a profile for $(4,2,4,...,4,4,2)$ where there are $2k$ (i.e., an even number
of) occurrences of $4$.\medskip

\textbf{Case 2.} There are $2k+1$ (i.e., odd number of) occurrences of $4$%
.\medskip

We follow the same steps as in case 2 of sequence I until for individual \#1
at profile $v$, for the options $\{x_{1},x_{2},\dots ,x_{2k+3}\}$, we change
the initial ranking of $x_{1}\succ x_{2}\succ \cdots \succ x_{2k+3}$
to:\medskip 
\[
\begin{tabular}{|l|}
\hline
$%
\color{red}%
x_{2k+2}$ \\ 
$%
\color{red}%
x_{2k+3}$ \\ 
$%
\color{green}%
x_{2k+1}$ \\ 
$x_{2k-1}$ \\ 
$x_{2k}$ \\ 
$\vdots $ \\ 
$x_{3}$ \\ 
$x_{4}$ \\ 
$x_{1}$ \\ 
$x_{2}$ \\ \hline
\end{tabular}%
\]

The options are split into $(k+2)$ levels, with the Borda score change from $%
v$ to $u$ marked as below:\medskip 
\begin{eqnarray*}
\{x_{2k+2},x_{2k+3}\} &:&(2k+1)\,more \\
\{x_{2k+1}\} &:&(2k-2)\,more \\
\{x_{2k-1},x_{2k}\} &:&(2k-5)\,more \\
&&\vdots \\
\{x_{3},x_{4}\} &:&(3-2k)\,more\,[i.e.,(2k-3)\,less] \\
\{x_{1},x_{2}\} &:&(-1-2k)\,more\,[i.e.,(2k+1)\,less]
\end{eqnarray*}

The Borda score difference between the black options based on the ranking at 
$u$ and the blue options based on the ranking at $v$ (therefore a single
level for the blue options) is:\medskip 
\begin{eqnarray*}
\{x_{2k+2},x_{2k+3}\} &:&(2k+1)+(2k+2)\,more\,[i.e.,(4k+3)\,more] \\
\{x_{2k+1}\} &:&(2k-2)+(2k+2)\,more\,[i.e.,4k\,more] \\
\{x_{2k-1},x_{2k}\} &:&(2k-5)+(2k+2)\,more\,[i.e.,(4k-3)\,more] \\
&&\vdots \\
\{x_{3},x_{4}\} &:&(3-2k)+(2k+2)\,more\,[i.e.,5\,more] \\
\{x_{1},x_{2}\} &:&(-1-2k)+(2k+2)\,more\,[i.e.,1\,more]
\end{eqnarray*}

For the options $\{x_{(2k+3)+1},x_{(2k+3)+2},\dots ,x_{(2k+3)+2k+1}\}$, we
make the exactly same change as case 2 of sequence I, so the options are
split into $(k+1)$ and the Borda score change from $v$ to $u$ is:\medskip 
\begin{eqnarray*}
\{%
\color{blue}%
x_{(2k+1)+2k},%
\color{blue}%
x_{(2k+1)+2k+1}%
\color{black}%
\} &:&(2k-1)\,more \\
\{%
\color{blue}%
x_{(2k+1)+2k-2},%
\color{blue}%
x_{(2k+1)+2k-1}%
\color{black}%
\} &:&(2k-5)\,more \\
&&\vdots \\
\{%
\color{blue}%
x_{(2k+1)+4},%
\color{blue}%
x_{(2k+1)+5}%
\color{black}%
\} &:&(7-2k)\,more\,[i.e.,(2k-7)\,less] \\
\{%
\color{blue}%
x_{(2k+1)+2},%
\color{blue}%
x_{(2k+1)+3}%
\color{black}%
\} &:&(3-2k)\,more\,[i.e.,(2k-3)\,less] \\
\{%
\color{blue}%
x_{(2k+1)+1}%
\color{black}%
\} &:&-2k\,more\,[i.e.,2k\,less]
\end{eqnarray*}

We then compare the above $(k+2)$-level group and the $(k+1)$-level
group.\medskip

For $k\geq 1$, we have $4k>2k-1$ (actually, we have $4k>2k+1$ also, if we
need to make the adjustments below), so the two levels consisting of $%
\{x_{2k+2},x_{2k+3}\}$ and $\{x_{2k+1}\}$ remain to be the top two levels
among the $(2k+3)$ levels. Similarly, the level $\{%
\color{blue}%
x_{(2k+1)+1}%
\color{black}%
\}$ remains to be the lowest one among the $(2k+3)$ levels. Suppose there
exists a level among the $(k+1)$-level group of blue options where the Borda
score change for the options at the level from $v$ to $u$ is equal to $1$
(if no such level exists, then we skip this step), we do the same
adjustments as illustrated at the end of case 2 of sequence I. The obtained
profile is for $(4,2,4,...,4,4,2)$ where there are $2k+1$ (i.e., an odd
number of) occurrences of $4$.\bigskip

\textbf{III. The sequence }$(2,4,...,4,4,2,4)$.\medskip

The inverse of the sequence $(2,4,...,4,4,2,4)$ is $(4,2,4,4,...,4,2)$,
which we have shown above is in the Borda range for $n=3$ (and thus $n\geq 3$%
). \ Taking the profile for $(4,2,4,4,...,4,2)$ and inverting everyone's
ranking yields the sequence $(2,4,...,4,4,2,4)$.\bigskip

\textbf{IV. The sequence }$(4,2,4,...,4,2,4)$.\medskip

\textbf{Case 1.} There are $2k$ (i.e., even number of) $4s$. Note that when $%
k=1$, the sequence is $(4,2,2,4)$, which is shown to be in the Borda range
in the Appendix.\ So we focus on $k>1$.\medskip

For the options, $x_{1},x_{2},\dots ,x_{2k+1}$ (black), we follow the same
steps as in sequence II, case 1, and accordingly, the Borda score difference
between the black options based on ranking at $u$ and the blue options based
on ranking at $v$ (therefore a single level for the blue options)
is:\medskip 
\begin{eqnarray*}
\{x_{2k},x_{2k+1}\} &:&(2k-1)+(2k+1)\,more\,[i.e.,4k\,more] \\
\{x_{2k-1}\} &:&(2k-4)+(2k+1)\,more\,[i.e.,(4k-3)\,more] \\
\{x_{2k-3},x_{2k-2}\} &:&(2k-7)+(2k+1)\,more\,[i.e.,(4k-6)\,more] \\
\{x_{2k-5},x_{2k-4}\} &:&(2k-11)+(2k+1)\,more\,[i.e.,(4k-10)\,more] \\
&&\vdots \\
\{x_{3},x_{4}\} &:&(5-2k)+(2k+1)\,more\,[i.e.,6\,more] \\
\{x_{1},x_{2}\} &:&(1-2k)+(2k+1)\,more\,[i.e.,2\,more]
\end{eqnarray*}

For the options $\{x_{(2k+1)+1},\dots ,x_{(2k+1)+2k+1}\}$, we follow the
same steps as case 1 of sequence I until the part where we change the
initial ranking of $x_{(2k+1)+1}\succ x_{(2k+1)+2}\succ \cdots \succ
x_{(2k+1)+2k+1}$ to: 
\[
\begin{tabular}{|l|}
\hline
$%
\color{blue}%
x_{(2k+1)+2k}$ \\ 
$%
\color{blue}%
x_{(2k+1)+2k+1}$ \\ 
$%
\color{blue}%
x_{(2k+1)+2k-2}$ \\ 
$%
\color{blue}%
x_{(2k+1)+2k-1}$ \\ 
$%
\color{blue}%
\vdots $ \\ 
$%
\color{blue}%
x_{(2k+1)+4}$ \\ 
$%
\color{blue}%
x_{(2k+1)+5}$ \\ 
$%
\color{green}%
x_{(2k+1)+3}$ \\ 
$%
\color{red}%
x_{(2k+1)+1}$ \\ 
$%
\color{red}%
x_{(2k+1)+2}$ \\ \hline
\end{tabular}%
\]

And therefore, the options are split into $(k+1)$ levels and\ their Borda
score change is:\medskip 
\begin{eqnarray*}
\{%
\color{blue}%
x_{(2k+1)+2k},%
\color{blue}%
x_{(2k+1)+2k+1}%
\color{black}%
\} &:&(2k-1)\,more \\
\{%
\color{blue}%
x_{(2k+1)+2k-2},%
\color{blue}%
x_{(2k+1)+2k-1}%
\color{black}%
\} &:&(2k-5)\,more \\
&&\vdots \\
\{%
\color{blue}%
x_{(2k+1)+4},%
\color{blue}%
x_{(2k+1)+5}%
\color{black}%
\} &:&(7-2k)\,more\,[i.e.,(2k-7)\,less] \\
\{%
\color{blue}%
x_{(2k+1)+3}%
\color{black}%
\} &:&(4-2k)\,more\,[i.e.,(2k-4)\,less] \\
\{%
\color{blue}%
x_{(2k+1)+1},%
\color{blue}%
x_{(2k+1)+1}%
\color{black}%
\} &:&(1-2k)\,more\,[i.e.,(2k-1)\,less]
\end{eqnarray*}

We compare the earlier $(k+1)$-level group and the above $(k+1)$-level
group.\medskip

Recall that we focus on $k>1$, so $(4k-3)>(2k-1)$, and therefore, the two
levels $\{x_{2k},x_{2k+1}\}$ and $\{x_{2k-1}\}$ are still the highest two
levels among the total $2(k+1)$ levels. Similarly, the two levels $\{%
\color{blue}%
x_{(2k+1)+3}%
\color{black}%
\}$ and $\{%
\color{blue}%
x_{(2k+1)+1},%
\color{blue}%
x_{(2k+1)+1}%
\color{black}%
\}$ remain to be the lowest levels among the total $2(k+1)$ levels. For the
remaining options, the Borda score difference for black options is even
while for blue the difference is odd, so they still remain to be $2k$
levels. Thus, the obtained profile is for $(4,2,4,...,4,2,4)$ where there
are $2k$ (i.e., even number of) $4s$.\medskip

\textbf{Case 2.} There are $2k+1$ (i.e., odd number of) occurrences of $4$%
.\medskip

For options $\{x_{1},x_{2},\dots ,x_{2k+3}\}$, we follow exactly the same
steps as for case 2 of sequence II, where the Borda score difference between
these options based on ranking at $u$ and the blue options based on ranking
at $v$ (therefore a single level for the blue options) is:\medskip 
\begin{eqnarray*}
\{x_{2k+2},x_{2k+3}\} &:&(2k+1)+(2k+2)\,more\,[i.e.,(4k+3)\,more] \\
\{x_{2k+1}\} &:&(2k-2)+(2k+2)\,more\,[i.e.,4k\,more] \\
\{x_{2k-1},x_{2k}\} &:&(2k-5)+(2k+2)\,more\,[i.e.,(4k-3)\,more] \\
&&\vdots \\
\{x_{3},x_{4}\} &:&(3-2k)+(2k+2)\,more\,[i.e.,5\,more] \\
\{x_{1},x_{2}\} &:&(-1-2k)+(2k+2)\,more\,[i.e.,1\,more]
\end{eqnarray*}

For options $\{x_{(2k+3)+1},x_{(2k+3)+2},\dots ,x_{(2k+3)+2k+1}\}$ here, we
follow the steps for options $\{x_{(2k+1)+1},\dots ,x_{(2k+1)+2k+1}\}$ in
the above sequence IV, case 1. And therefore, the options are split into $%
(k+1)$ levels and\ their Borda score change is:\medskip 
\begin{eqnarray*}
\{%
\color{blue}%
x_{(2k+3)+2k},%
\color{blue}%
x_{(2k+3)+2k+1}%
\color{black}%
\} &:&(2k-1)\,more \\
\{%
\color{blue}%
x_{(2k+3)+2k-2},%
\color{blue}%
x_{(2k+3)+2k-1}%
\color{black}%
\} &:&(2k-5)\,more \\
&&\vdots \\
\{%
\color{blue}%
x_{(2k+3)+4},%
\color{blue}%
x_{((2k+3)+5}%
\color{black}%
\} &:&(7-2k)\,more\,[i.e.,(2k-7)\,less] \\
\{%
\color{blue}%
x_{(2k+3)+3}%
\color{black}%
\} &:&(4-2k)\,more\,[i.e.,(2k-4)\,less] \\
\{%
\color{blue}%
x_{(2k+3)+1},%
\color{blue}%
x_{(2k+3)+1}%
\color{black}%
\} &:&(1-2k)\,more\,[i.e.,(2k-1)\,less]
\end{eqnarray*}

Again we compare the two groups.\medskip

We focus here on $k\geq 2$. The cases where $k=0,1$ will be treated
separately in the Appendix.\medskip

For $k\geq 2$, we have $4-2k<1$ and $4k>2k-1$ (actually we even have $%
4k>2k+1 $ if we need to make the adjustments below). So $\{x_{2k+2},x_{2k+3}%
\}$ and $\{x_{2k+1}\}$ remain to be the highest levels among the $(2k+3)$
levels. And $\{%
\color{blue}%
x_{(2k+3)+3}%
\color{black}%
\}$ and $\{%
\color{blue}%
x_{(2k+3)+1},%
\color{blue}%
x_{(2k+3)+1}%
\color{black}%
\}$ remain to be the lowest levels among the $(2k+3)$ levels. If there
exists a level among the $(k+1)$-level group of blue options where the Borda
score change for the options at the level from $v$ to $u$ is equal to $1$
(if no such level exists, then we skip this step), we do the same
adjustments as illustrated at the end of sequence I, case 2. The obtained
profile is for $(4,2,4,...,4,2,4)$ where there are $2k+1$ (i.e., odd number
of) $4s$.\medskip

\textbf{4. \ Remark.}\medskip

A modification of the analysis above allows us to prove an extended version
of the lemma:

\qquad For positive integer $t$, each of the following (patterns of) orders
(with exactly two levels equal to $2$) is in the Borda range for $n=3$ (and
so for all odd $n\geq 3$):\medskip

\qquad \qquad \qquad \qquad $(2,4t,4t,...,4t,4t,2)$

\qquad \qquad \qquad \qquad $(4t,2,4t,4t,...,4t,4t,2)$

\qquad \qquad \qquad \qquad $(2,4t,4t,...,4t,4t,2,4t)$

\qquad \qquad \qquad \qquad $(4t,2,4t,4t,...,4t,4t,2,4t)$\medskip

\textbf{Appendix}\medskip

We complete the proof of the Lemma by showing that the following sequences
(4,2,4,2), (4,2,2,4), (4,2,2), (2,2,4) and (4,2,4,2,4) are in the Borda
range for $n=3$.\medskip

The following profile shows that $(4,2,4,2)$ is in the Borda range for $n=3$%
, where the four levels consist of $\{x_{4},x_{5},x_{10},x_{11}\}$, $%
\{x_{1},x_{7}\}$, $\{x_{2},x_{3},x_{8},x_{9}\}$ and $\{x_{6},x_{12}\}$%
.\medskip 
\[
\begin{tabular}{|lll|}
\hline
$1$ & $2$ & $3$ \\ \hline
$x_{1}$ & \multicolumn{1}{|l}{$x_{10}$} & \multicolumn{1}{|l|}{$x_{11}$} \\ 
$x_{4}$ & \multicolumn{1}{|l}{$x_{8}$} & \multicolumn{1}{|l|}{$x_{9}$} \\ 
$x_{5}$ & \multicolumn{1}{|l}{$x_{12}$} & \multicolumn{1}{|l|}{$x_{7}$} \\ 
$x_{2}$ & \multicolumn{1}{|l}{$x_{4}$} & \multicolumn{1}{|l|}{$x_{5}$} \\ 
$x_{3}$ & \multicolumn{1}{|l}{$x_{2}$} & \multicolumn{1}{|l|}{$x_{3}$} \\ 
$x_{6}$ & \multicolumn{1}{|l}{$x_{6}$} & \multicolumn{1}{|l|}{$x_{1}$} \\ 
$x_{7}$ & \multicolumn{1}{|l}{$x_{11}$} & \multicolumn{1}{|l|}{$x_{12}$} \\ 
$x_{10}$ & \multicolumn{1}{|l}{$x_{9}$} & \multicolumn{1}{|l|}{$x_{10}$} \\ 
$x_{11}$ & \multicolumn{1}{|l}{$x_{7}$} & \multicolumn{1}{|l|}{$x_{8}$} \\ 
$x_{8}$ & \multicolumn{1}{|l}{$x_{5}$} & \multicolumn{1}{|l|}{$x_{6}$} \\ 
$x_{9}$ & \multicolumn{1}{|l}{$x_{3}$} & \multicolumn{1}{|l|}{$x_{4}$} \\ 
$x_{12}$ & \multicolumn{1}{|l}{$x_{1}$} & \multicolumn{1}{|l|}{$x_{2}$} \\ 
\hline
\end{tabular}%
\]%
\bigskip

The next profile shows that $(4,2,2,4)$ is in the Borda range for $n=3$,
where the four levels consist of $\{x_{2},x_{3},x_{8},x_{9}\}$, $%
\{x_{6},x_{12}\}$, $\{x_{1},x_{7}\}$ and $\{x_{4},x_{5},x_{10},x_{11}\}$%
.\medskip 
\[
\begin{tabular}{|lll|}
\hline
$1$ & $2$ & $3$ \\ \hline
$x_{2}$ & \multicolumn{1}{|l}{$x_{8}$} & \multicolumn{1}{|l|}{$x_{9}$} \\ 
$x_{3}$ & \multicolumn{1}{|l}{$x_{12}$} & \multicolumn{1}{|l|}{$x_{7}$} \\ 
$x_{1}$ & \multicolumn{1}{|l}{$x_{10}$} & \multicolumn{1}{|l|}{$x_{11}$} \\ 
$x_{6}$ & \multicolumn{1}{|l}{$x_{2}$} & \multicolumn{1}{|l|}{$x_{3}$} \\ 
$x_{4}$ & \multicolumn{1}{|l}{$x_{6}$} & \multicolumn{1}{|l|}{$x_{1}$} \\ 
$x_{5}$ & \multicolumn{1}{|l}{$x_{4}$} & \multicolumn{1}{|l|}{$x_{5}$} \\ 
$x_{8}$ & \multicolumn{1}{|l}{$x_{11}$} & \multicolumn{1}{|l|}{$x_{12}$} \\ 
$x_{9}$ & \multicolumn{1}{|l}{$x_{9}$} & \multicolumn{1}{|l|}{$x_{10}$} \\ 
$x_{7}$ & \multicolumn{1}{|l}{$x_{7}$} & \multicolumn{1}{|l|}{$x_{8}$} \\ 
$x_{12}$ & \multicolumn{1}{|l}{$x_{5}$} & \multicolumn{1}{|l|}{$x_{6}$} \\ 
$x_{10}$ & \multicolumn{1}{|l}{$x_{3}$} & \multicolumn{1}{|l|}{$x_{4}$} \\ 
$x_{11}$ & \multicolumn{1}{|l}{$x_{1}$} & \multicolumn{1}{|l|}{$x_{2}$} \\ 
\hline
\end{tabular}%
\]%
\bigskip

We move on to $(4,2,2)$ and $(2,2,4)$. Note that we only need to construct a
profile for one of them and inverting everyone's ranking yields the other.
The following profile shows that $(4,2,2)$ is in the Borda range for $n=3$,
where the three levels consist of $\{x_{2},x_{3},x_{6},x_{7}\}$, $%
\{x_{4},x_{8}\}$ and $\{x_{1},x_{5}\}$.\medskip 
\[
\begin{tabular}{|lll|}
\hline
$1$ & $2$ & $3$ \\ \hline
$x_{2}$ & \multicolumn{1}{|l|}{$x_{7}$} & \multicolumn{1}{|l|}{$x_{6}$} \\ 
$x_{3}$ & \multicolumn{1}{|l|}{$x_{5}$} & \multicolumn{1}{|l|}{$x_{8}$} \\ 
$x_{4}$ & \multicolumn{1}{|l|}{$x_{3}$} & \multicolumn{1}{|l|}{$x_{2}$} \\ 
$x_{1}$ & \multicolumn{1}{|l|}{$x_{1}$} & \multicolumn{1}{|l|}{$x_{4}$} \\ 
$x_{6}$ & \multicolumn{1}{|l}{$x_{8}$} & \multicolumn{1}{|l|}{$x_{7}$} \\ 
$x_{7}$ & \multicolumn{1}{|l}{$x_{6}$} & \multicolumn{1}{|l|}{$x_{5}$} \\ 
$x_{8}$ & \multicolumn{1}{|l}{$x_{4}$} & \multicolumn{1}{|l|}{$x_{3}$} \\ 
$x_{5}$ & \multicolumn{1}{|l}{$x_{2}$} & \multicolumn{1}{|l|}{$x_{1}$} \\ 
\hline
\end{tabular}%
\]%
\bigskip

The last profile shows that $(4,2,4,2,4)$ is also in the Borda range for $%
n=3 $, where the five levels consist of $\{x_{4},x_{5},x_{12},x_{13}\}$, $%
\{x_{3},x_{11}\}$, $\{x_{7},x_{8},x_{15},x_{16}\}$, $\{x_{6},x_{14}\}$ and $%
\{x_{1},x_{2},x_{9},x_{10}\}$.\medskip 
\[
\begin{tabular}{|lll|}
\hline
$1$ & $2$ & $3$ \\ \hline
\multicolumn{1}{|l|}{$x_{4}$} & \multicolumn{1}{|l}{$x_{12}$} & 
\multicolumn{1}{|l|}{$x_{13}$} \\ 
\multicolumn{1}{|l|}{$x_{5}$} & \multicolumn{1}{|l}{$x_{10}$} & 
\multicolumn{1}{|l|}{$x_{11}$} \\ 
\multicolumn{1}{|l|}{$x_{3}$} & \multicolumn{1}{|l}{$x_{16}$} & 
\multicolumn{1}{|l|}{$x_{9}$} \\ 
\multicolumn{1}{|l|}{$x_{7}$} & \multicolumn{1}{|l}{$x_{14}$} & 
\multicolumn{1}{|l|}{$x_{15}$} \\ 
\multicolumn{1}{|l|}{$x_{8}$} & \multicolumn{1}{|l}{$x_{4}$} & 
\multicolumn{1}{|l|}{$x_{5}$} \\ 
$x_{6}$ & \multicolumn{1}{|l}{$x_{2}$} & \multicolumn{1}{|l|}{$x_{3}$} \\ 
$x_{1}$ & \multicolumn{1}{|l}{$x_{8}$} & \multicolumn{1}{|l|}{$x_{1}$} \\ 
$x_{2}$ & \multicolumn{1}{|l}{$x_{6}$} & \multicolumn{1}{|l|}{$x_{7}$} \\ 
$x_{12}$ & \multicolumn{1}{|l}{$x_{15}$} & \multicolumn{1}{|l|}{$x_{16}$} \\ 
$x_{13}$ & \multicolumn{1}{|l}{$x_{13}$} & \multicolumn{1}{|l|}{$x_{14}$} \\ 
$x_{11}$ & \multicolumn{1}{|l}{$x_{11}$} & \multicolumn{1}{|l|}{$x_{12}$} \\ 
$x_{15}$ & \multicolumn{1}{|l}{$x_{9}$} & \multicolumn{1}{|l|}{$x_{10}$} \\ 
$x_{16}$ & \multicolumn{1}{|l}{$x_{7}$} & \multicolumn{1}{|l|}{$x_{8}$} \\ 
$x_{14}$ & \multicolumn{1}{|l}{$x_{5}$} & \multicolumn{1}{|l|}{$x_{6}$} \\ 
$x_{9}$ & \multicolumn{1}{|l}{$x_{3}$} & \multicolumn{1}{|l|}{$x_{4}$} \\ 
$x_{10}$ & \multicolumn{1}{|l}{$x_{1}$} & \multicolumn{1}{|l|}{$x_{2}$} \\ 
\hline
\end{tabular}%
\]%
\bigskip \bigskip

\textbf{References}\medskip

[1] \ Kelly, JS and S Qi (2015) "The Construction of Orderings by Borda's
Rule".\medskip 

\end{document}